\documentclass[12pt,twoside]{article}
\usepackage{graphicx}
\usepackage{amsfonts,amsbsy,amsthm,amssymb,amsmath,tikz}
\allowdisplaybreaks
\usepackage{cite}
\usepackage{graphics}
\def\bpsp{\begin{pspicture}}
\def\epsp{\end{pspicture}}

\newtheorem{theorem}{Theorem}[section]
\newtheorem{remark}[theorem]{Remark}
\newtheorem{example}[theorem]{Example}
\newtheorem{lemma}[theorem]{Lemma}
\newtheorem{corollary}[theorem]{Corollary}
\newtheorem{definition}[theorem]{Definition}
\newtheorem{proposition}[theorem]{Proposition}

\newtheorem{note}{Note}
\newtheorem{fact}{Fact}
\newtheorem{case}{Case}
\newtheorem{conjecture}{Conjecture}
\newtheorem{question}{Question}

\newcommand{\bea}{\begin{eqnarray}}
\newcommand{\eea}{\end{eqnarray}}
\newcommand{\beq}{\begin{eqnarray*}}
\newcommand{\eeq}{\end{eqnarray*}}

\def\m4{\mbox{\rm ~(mod $4$)}}

\def \bd{\begin{definition}}
\def \ed{\end{definition}}
\def \bqu{\begin{question}}
\def \equ{\end{question}}
\def \bcc{\begin{conjecture}}
\def \ecc{\end{conjecture}}
\def \bt{\begin{theorem}}
\def \et{\end{theorem}}
\def \bl{\begin{lemma}}
\def \el{\end{lemma}}
\def \bc{\begin{corollary}}
\def \ec{\end{corollary}}
\def \be{\begin{equation}}
\def \ee{\end{equation}}
\def \ben{\begin{enumerate}}
\def \een{\end{enumerate}}
\def \ba{\begin{array}}
\def \ea{\end{array}}
\def \bp{\begin{proposition}}
\def \ep{\end{proposition}}
\def \bx{\begin{example}}
\def \ex{\end{example}}
\def \br{\begin{remark}}
\def \er{\end{remark}}
\def \bdsc{\begin{description}}
\def \edsc{\end{description}}

\def \bn{\begin{case}}
\def \en{\end{case}}
\def \bnt{\begin{note}}
\def \ent{\end{note}}
\def\1{1\!\!1}

\def\mm2{\mbox{\rm ~(mod $2$)}}
\def\m4{\mbox{\rm ~(mod $4$)}}

\def\qed{\nolinebreak\hfill\rule{.2cm}{.2cm}\par\addvspace{.5cm}}

\def\m{\mu}

\def\1{\textbf{1}}
\def\0{\textbf{0}}

\begin{document}
\title{On graphs with distance Laplacian eigenvalues of multiplicity $n-4$}
\author{ Saleem Khan$ ^{a} $, S. Pirzada$ ^{b} $\\
$^{a,b}${\em Department of Mathematics, University of Kashmir, Srinagar, Kashmir, India}\\
$ ^{a} $\texttt{khansaleem1727@gmail.com}; $^{b}$\texttt{pirzadasd@kashmiruniversity.ac.in}
}
\date{}

\pagestyle{myheadings} \markboth{Khan,Pirzada}{On graphs with distance Laplacian eigenvalues of multiplicity $n-4$}
\maketitle
\vskip 5mm
\noindent{\footnotesize \bf Abstract.} Let $G$ be a connected simple graph with $n$ vertices. The distance Laplacian matrix $D^{L}(G)$ is defined as $D^L(G)=Diag(Tr)-D(G)$, where $Diag(Tr)$ is the diagonal matrix of vertex transmissions and $D(G)$ is the distance matrix of $G$. The eigenvalues of $D^{L}(G)$ are the distance Laplacian eigenvalues of $G$ and are denoted by $\partial_{1}^{L}(G)\geq \partial_{2}^{L}(G)\geq \dots \geq \partial_{n}^{L}(G)$. The largest eigenvalue $\partial_{1}^{L}(G)$ is called the distance Laplacian spectral radius. Lu et al. (2017), Fernandes et al. (2018) and Ma et al. (2018) completely characterized the graphs having some distance Laplacian eigenvalue of multiplicity $n-3$. In this paper, we characterize the graphs having distance Laplacian spectral radius of multiplicity $n-4$ together with one of the distance Laplacian eigenvalue as $n$ of multiplicity either 3 or 2. Further, we completely determine the graphs for which the distance Laplacian eigenvalue $n$ is of multiplicity $n-4$.

\vskip 3mm



\section{Introduction}\label{sec1}

Throughout this paper, we consider simple and connected graphs. A simple connected graph  $G=(V,E)$ consists of the vertex set $V(G)=\{v_{1},v_{2},\ldots,v_{n}\}$ and the edge set  $E(G)$. The \textit{order} and \textit{size} of $G$ are $|V(G)|=n$ and  $|E(G)|=m$, respectively. The \textit{degree} of a vertex $v,$ denoted by $d_{G}(v)$ (we simply write by $d_v$) is the number of edges incident on the vertex $v$.  For other standard definitions, we refer to \cite{A1,sp}. The adjacency matrix $A=(a_{ij})$ of $G$ is an $n\times n$ matrix whose $(i,j)$-entry is equal to 1, if $v_i$ is adjacent to $v_j$ and equal to $ 0 $, otherwise. Let $Deg(G)=\text{diag}(d_{v_1}(G), d_{v_2}(G), \dots, d_{v_n}(G))$ be the diagonal matrix of vertex degrees $d_{v_i}(G)$, $i=1,2,\dots,n$. The positive semi-definite matrix $L(G)=Deg(G)-A(G)$ is the Laplacian matrix of $G$. The eigenvalues of $L(G)$ are called the Laplacian eigenvalues of $G$. The Laplacian eigenvalues are denoted by $\mu_1 (G) ,\mu_2 (G),\dots,\mu_n (G)$ and are ordered as $\mu_1 (G) \geq \mu_2 (G)\geq \dots \geq \mu_n (G)$. The sequence of the Laplacian eigenvalues is called the Laplacian spectrum (briefly $L$-spectrum) of $G$. In $G$, the \textit{distance} between two vertices $u,v\in V(G),$ denoted by $d_{uv}=d(u,v)$, is defined as the length of a shortest path between $u$ and $v$. The diameter of $G$ denoted by $diam(G)$ is $\max_{u,v\in G}d(u,v)$, that is, the length of a longest path among the distance between every two vertices of $G$.  The \textit{distance matrix} of $G$ is defined as $D(G)=(d_{uv})_{u,v\in V(G)}$. The \textit{transmission} $Tr_{G}(v)$ of a vertex $v$ is the sum of the distances from $v$ to all other vertices in $G$, that is, $Tr_{G}(v)=\sum\limits_{u\in V(G)}d_{uv}.$ For any vertex $v_i \in V(G)$, the transmission $Tr_G(v_i)=Tr_i$ is also called the \textit{transmission degree}.

Let $Diag(Tr)=Diag (Tr_1,Tr_2,\ldots,Tr_n) $ be the diagonal matrix of vertex transmissions of $G$. Aouchiche and Hansen \cite{A2} defined the \textit{distance Laplacian matrix} of $G$ as $D^L(G)=Diag(Tr)-D(G)$ (or simply $D^{L}$). The eigenvalues of $D^{L}$ are called the distance Laplacian eigenvalues of $G$. Clearly, $ D^L(G) $ is a real symmetric positive semi-definite matrix so that its eigenvalues can be ordered as $ \partial^L_{1}(G)\geq \partial^L_{2}(G)\geq \dots\geq \partial^L _{n-1}(G)>\partial^L _{n}(G)=0 $.  We write $\partial^L_{i}$ in place of $\partial^L_{i}(G)$ if the graph $G$ is clear from the context. If $G$ has $k$ distinct distance Laplacian eigenvalues say $\partial^L_{1}(G),\partial^L_{2}(G),\dots,\partial^L_{k}(G)$ with corresponding multiplicities as $n_1 ,n_2 ,\dots,n_k$,  we write the  distance Laplacian spectrum ( briefly $D^{L}$-spectrum) of $G$  as $\Big({\partial^{L}_{1}}^{(n_1)},{\partial^{L}_{2}}^{(n_2)},\dots,{\partial^{L}_{k}}^{(n_k)}\Big)$. The largest eigenvalue $\partial^L _{1}(G) $ is called the distance Laplacian spectral radius of $G$. We denote the multiplicity of the distance Laplacian eigenvalue $\partial^L_{i}(G)$ by $m(\partial^L_{i}(G))$. More recent work on distance Laplacian eigenvalues can be seen in \cite{ps}.\\
\indent As usual,  $K_n$, $C_n$, $P_n$ and $S_n$ are respectively, the complete graph, the cycle, the path and the star all on $n$ vertices. A clique of a graph $G$ is an induced subgraph of $G$ that is complete. A $kite$ $Ki_{n,\omega}$ is the graph obtained from a clique $K_\omega$ and a path$P_{n-\omega}$ by adding an edge between an endpoint of the path and a vertex from the clique. $SK_{n,\alpha}$ denotes the complete split graph, that is, the complement of the disjoint union of a clique $K_\alpha$ and $n-\alpha$ isolated vertices. A complete mutipartite graph is denoted by $K_{t_1 ,t_2 ,\dots,t_l}$, where $l$ is the number of partite classes and $t_1 +t_2 +\dots+t_l =n$. Throughout, we assume that $t_1 \geq t_2 \geq \dots \geq t_l$. If $l=2$, it is a complete bipartite graph. Usually we will choose the complete bipartite graph $K_{n-r,~r}$, for $1\leq r\leq \frac{n}{2}$.  If $G$ is a non complete graph on $n\geq 2$ vertices, then $G+e$ is the graph obtained from $G$ by adding an edge $e$ between any two non-adjacent vertices. Further, if $f$ be an edge of $G$, then $G-f$ is the graph obtained from $G$ by deleting the edge $f$. The join of two graphs $G_1$ and $G_2$, denoted by $G_1\vee G_2$, is a graph obtained from $G_1$ and $G_2$ by joining each vertex of $G_1$ to all vertices of $G_2$. \\
\indent Fernandes et al. \cite{A3} and Lu et al. \cite{A6} determined the graphs having distance Laplacian spectral radius of multiplicity $n-3$. Further, the investigation of the graphs having some distance Laplacian eigenvalue of multiplicity $n-3$ was done by   Ma et al. in \cite{A7}. \\
\indent  The rest of the paper is organized as follows. In Section $2$, we state some preliminary results, which will be used to prove our main results. In Section 3, we characterize the graphs having distance Laplacian spectral radius of multiplicity $n-4$ together with one of the distance Laplacian eigenvalue as $n$ of multiplicity either 3 or 2. In Section 4, we completely determine the graphs for which the distance Laplacian eigenvalue $n$ is of multiplicity $n-4$.

\section{Preliminaries }

\begin{lemma}\label{L1}\emph{\cite{A1}} Let $G$ be a graph on n vertices with Laplacian eigenvalues $\mu_1 (G) \geq \mu_2 (G)\geq \dots \geq \mu_n (G)=0$. Then the Laplacian eigenvalues of $\overline{G}$ are given by $\mu_i (\overline{G})=n-\mu_{n-i} (G)$ for $i=1,\dots,n-1$ and $\mu_n (\overline{G})=0$.
\end{lemma}

\begin{lemma}\label{L2}\emph{\cite{A2}} Let $G$ be a connected graph on $n$ vertices with $diam(G)\leq 2$. Let $\mu_1 (G) \geq \mu_2 (G)\geq \dots \geq \mu_n (G)=0$ be the Laplacian eigenvalues of $G$. Then the distance Laplacian eigenvalues of $G$ is  $2n-\mu_{n-1} (G) \geq 2n- \mu_{n-2} (G)\geq \dots \geq 2n-\mu_1 (G)>\partial^L_n (G)=0$. Moreover, for every $i\in \{1,2,\dots,n-1\}$ the eigenspaces corresponding to $\mu_i (G)$ and $2n-\mu_i (G)$ are the same.
\end{lemma}

\begin{lemma}\label{L3}\emph{\cite{A2}} Let $G$ be a connected graph on $n$ vertices. Then $\partial^L_{n-1}(G)=n$ if and only if $\overline {G}$ is disconnected. Furthermore, the multiplicity of $n$ as a distance Laplacian eigenvalue is one less than the number of connected components of $\overline {G}$.
\end{lemma}

\begin{lemma}\label{L4} \emph{\cite{A8}}  Let $t_{1},t_{2},\dots,t_{k}$ and n be integers such that $t_{1}+t_{2}+\dots+t_{k}=n$ and $t_{i}\geq 1$ for $i=1,2,\dots,k$. Let $p=|\{i:t_{i}\geq 2\}|$. The distance Laplacian spectrum of the complete $k-partite$ graph $K_{t_{1},t_{2},\dots,t_{k}}$ is $ \Big({(n+t_{1})}^{(t_{1}-1)},\dots,{(n+t_{p})}^{(t_{p}-1)},n^{(k-1)},0\Big)$.
\end{lemma}

\begin{lemma}\label{L5} \emph{\cite{A3}} Let $G$ be a connected graph of order $n\geq 4$. Then, $G\cong K_n-e $ if and only if the $L-spectrum$ of $G$ is $\Big(\mu_1 ^{(n-2)},\mu_2 ,0\Big)$, with $\mu_1 > \mu_2 >0$.
\end{lemma}
 Let  $v\in V(G)$. By $N(v)$ we mean the set of all vertices which are adjacent to $v$ in $G$.
\begin{lemma}\label{L6} \emph{\cite{A5}} Let $G$ be a graph with $n$ vertices. If $K=\{v_1 ,v_2 ,\dots,v_p\}$ is a clique of $G$ such that $N(v_i)-K=N(v_j)-K$ for all $i,j\in \{1,2,\dots,p\}$, then $\partial=Tr(v_i)=Tr(v_j)$ for all $i,j\in \{1,2,\dots,p\}$ and $\partial +1$ is an eigenvalue of $D^L (G)$ with multiplicity at least $p-1$.
\end{lemma}

\begin{lemma}\label{L7} \emph{\cite{A5}} Let $G$ be a graph with $n$ vertices. If $K=\{v_1 ,v_2 ,\dots,v_p\}$ is an independent set of $G$ such that $N(v_i)=N(v_j)$ for all $i,j\in \{1,2,\dots,p\}$, then $\partial=Tr(v_i)=Tr(v_j)$ for all $i,j\in \{1,2,\dots,p\}$ and $\partial +2$ is an eigenvalue of $D^L (G)$ with multiplicity at least $p-1$.
\end{lemma}

\section{Multiplicity of distance Laplacian spectral radius}\label{sec2}

The following fact will be used frequently in the sequel.

\begin{fact}\label{fact}
A complete graph and a complete graph minus an edge are determined by their $L$-spectrum which is given by $\Big ((n)^{n-1},0\Big)$ and $\Big ((n)^{n-2},n-2,0\Big)$, respectively.
\end{fact}

Given a connected graph $G$ with order $n\geq 5$, we observe that one of the following possibilities can occur.\\
(a) $ m(\partial^L_1 (G))=n-4$ and $ \partial^L_{n-1} (G)=n$ with multiplicity 3,\\
(b) $ m(\partial^L_1 (G))=n-4$ and  $ \partial^L_{n-1} (G)=n$ with multiplicity 2,\\
(c) $ m(\partial^L_1 (G))=n-4$ and  $ \partial^L_{n-1} (G)=n$ with multiplicity 1,\\
(d) $ m(\partial^L_1 (G))=n-4$ and  $ \partial^L_{n-1} (G)\neq n$.\\

In the following theorem,  we address Cases (a) and (b), that is, we determine those graphs having the distance Laplacian spectral radius of multiplicity $n-4$ together with one of the distance Laplacian eigenvalue $n$ with multiplicity 3 or 2.

\begin{theorem}\label{T1}
Let $G$ be a connected graph with order $n\geq 6 $. Then\\
(a)  $ m(\partial^L_1 (G))=n-4$ and  $ \partial^L_{n-1} (G)=n$ with multiplicity 3 if and only if $G\cong K_{\frac{n}{4},\frac{n}{4},\frac{n}{4},\frac{n}{4}}$ if $n\equiv 0(mod~4)$, or $G\cong K_{\frac{n-1}{3},\frac{n-1}{3},\frac{n-1}{3},1}$ if $n\equiv 1(mod~3)$, or $G\cong K_{\frac{n-2}{2},\frac{n-2}{2},1,1}$ if $ n\equiv 2(mod~2)$ or $G\cong SK_{n,n-3}$.\\
(b)  $ m(\partial^L_1 (G))=n-4$ and  $ \partial^L_{n-1} (G)=n$ with multiplicity 2 if and only if $G\cong SK_{n,n-2}+e$ or $G\cong K_{n,n-3}-e $ or $G\cong K_{p,p,1}+e$; $p=\frac{n-1}{•2}\geq 3$ or $G\cong K_{p,p,2}$; $p=\frac{n-2}{•2}\geq 3$  or $G\cong K_{p,p,p}+e$; $p=\frac{n}{•3}\geq 3$.
\end{theorem}
\noindent{\bf Proof.} {\bf (a).} Using Lemma \ref{L3}, we note that $\overline{G}$ is disconnected having 4 components and $diam(G)=2$. Applying Lemmas \ref{L1} and \ref{L2}, the $L$-spectrum of  $\overline{G}$ is $\Big (({\partial^L_1 (G)-n})^{(n-4)},0,0,0,0\Big )$ so that every component of   $\overline{G}$ is either an isolated vertex or complete graphs with same order. Thus  $\overline{G}$ contains less or equal to three isolated vertices, that is,  $\overline{G}\cong K_{\frac{n}{4}}\cup K_{\frac{n}{4}}\cup K_{\frac{n}{4}}\cup K_{\frac{n}{4}}$ if $n\equiv 0(mod~4)$, or $\overline{G}\cong K_{\frac{n-1}{3}} \cup K_{\frac{n-1}{3}} \cup K_{ \frac{n-1}{3}}\cup K_1$ if $n-1\equiv 0(mod~3)$, or $\overline{G}\cong K_{\frac{n-2}{2}}\cup K_{\frac{n-2}{2}}\cup K_1 \cup K_1$ if $n-2\equiv 0(mod~2)$ or $\overline{G}\cong K_{n-3}\cup K_1 \cup K_1 \cup K_1$. Hence, $G\cong K_{\frac{n}{4},\frac{n}{4},\frac{n}{4},\frac{n}{4}}$, or $G\cong K_{\frac{n-1}{3},\frac{n-1}{3},\frac{n-1}{3},1}$, or $G\cong K_{\frac{n-2}{2},\frac{n-2}{2},1,1}$ according as $n\equiv 0(mod~4)$, or $n-1\equiv 0(mod~3)$, or $n-2\equiv 0(mod~2)$, or simply $G\cong SK_{n,n-3}$.\\
\indent Conversely, by the help of Lemma \ref{L4}, it is easy to see that the $D^L$-spectrum of $G\cong K_{\frac{n}{4},\frac{n}{4},\frac{n}{4},\frac{n}{4}}$,  $G\cong K_{\frac{n-1}{•3},\frac{n-1}{3},\frac{n-1}{3},1}$, $G\cong K_{\frac{n-2}{2},\frac{n-2}{2},1,1}$ and $G\cong SK_{n,n-3}$ are $\Big(\frac{5n}{4}^{(n-4)},n^{(3)},0\Big)$, $\Big(\frac{4n-1}{3}^{(n-4)},n^{(3)},0\Big)$, $\Big(\frac{3n-2}{2}^{(n-4)},n^{(3)},0\Big)$ and $\Big((2n-3)^{(n-4)},n^{(3)},0\Big)$, respectively.\\
\noindent{\bf (b).} For the graph $G$, let $n$ be a distance Laplacian eigenvalue with multiplicity 2. Using Lemma \ref{L3}, we see that $\overline{G}$ has three components, say $Q,~R$ and $S$, that is, $\overline{G}\cong Q\cup R\cup S$. This also shows that $diam(G)=2$. Assume that $|Q|\geq |R|\geq |S|$. By application of Lemmas \ref{L1} and \ref{L2}, we observe that the $L$-spectrum of $\overline{G}$ is $\Big(({\partial^L_1 (G)-n})^{(n-4)},\partial^L_{n-3} (G)-n,0,0,0\Big)$. We have the following possibilities.\\
\textbf{Case 1.} Let $|R|=|S|=1$. Then $L$-spectrum of $Q$ is $\Big({(\partial^L_1 (G)-n)}^{(n-4)},\partial^L_{n-3} (G)-n,0\Big)$. So, by Lemma \ref{L5}, $Q$ is isomorphic to $K_{n-4} \vee \overline{K_2}$. Therefore, $\overline{G}\cong (K_{n-4} \vee \overline{K_2})\cup K_1 \cup K_1$ which shows that $G\cong SK_{n,n-2}+e$.\\
\textbf{Case 2.} Let $|R|=2$, $|S|=1$. In order to find the $L$-spectrum of $\overline{G}$, we have to consider the following two subcases.\\
\textbf{Subcase 2.1.} Let the $L$-spectrum of $R$ be $\Big(\partial^L_{n-3} (G)-n,0\Big)$. Therefore, the $L$-spectrum of $Q$ is $\Big({(\partial^L_1 (G)-n)}^{(n-4)},0\Big)$. Using Fact \ref{fact}, from the above argument, it follows that $\partial^L_{n-3} (G)-n=2$ and $\partial^L_1 (G)-n=n-3$. Thus, $R\cong K_2$ and $Q\cong K_{n-3}$. From this, we obtain $\overline{G}\cong K_{n-3}\cup K_2 \cup K_1$ or $G\cong K_{n,n-3}-e $\\
 \textbf{Subcase 2.2.}  Let $L$-spectrum of $R$ be $\Big(\partial^L_1 (G)-n,0\Big)$. Therefore, the $L$-spectrum of $Q$ is $\Big({\partial^L_1 (G)-n}^{(n-5)},\partial^L_{n-3} (G)-n,0\Big)$. Thus, $R\cong K_2$. Using Lemma \ref{L5}, we have $Q\cong K_{n-5} \vee \overline{K_2}$. This shows that $\partial^L_1 (G)-n=2$. Also, $\partial^L_1 (G)-n=n-3$. Combining these two, we get $n=5$. This further implies that $\partial^L_1 (G)-n=\partial^L_{n-3} (G)-n$ or
$\partial^L_1 (G)=\partial^L_{n-3} (G)$, a contradiction.\\
\textbf{Case 3.} Let $|S|=1$, $|R|=p\geq 3$. To find the $L$-spectrum of $\overline{G}$, we see that either $L$-spectrum of $R$ is $\Big(({\partial^L_1 (G)-n})^{(p-1)},0\Big)$ and $L$-spectrum of $Q$ is $\Big(({\partial^L_1 (G)-n})^{(n-p-3)},\partial^L_{n-3} (G)-n,0\Big)$, or $L$-spectrum of $R$ is  $\Big(({\partial^L_1 (G)-n})^{(p-2)},\partial^L_{n-3} (G)-n,0\Big)$ and $L$-spectrum of $Q$ is $\Big(({\partial^L_1 (G)-n})^{(n-p-2)},0\Big)$. Using Fact \ref{fact}, we get $n=2p+1$, $\overline{G}\cong K_P \cup K_{p-2}\vee \overline{K_2}\cup K_1$ in both the cases. Hence $G\cong K_{p,p,1}+e$, $p=\frac{n-1}{•2}\geq 3$.   \\
\textbf{Case 4.}  Let $|S|= 2$, $|R|=p\geq 2$. We have to consider the following subcases.\\
\textbf{Subcase 4.1} Let $L$-spectrum of $S$ be $\Big({\partial^L_1 (G)-n},0\Big)$. Then we see that the $L$-spectrum of $R$ can be $\Big(({\partial^L_1 (G)-n})^{(p-1)},0\Big)$ or $\Big(({\partial^L_1 (G)-n})^{(n-p-4)},\partial^L_{n-3} (G)-n,0\Big)$. Again, using Fact \ref{fact}, we get  $\partial^L_1 (G)= \partial^L_{n-3} (G)$ in both the cases, which is a contradiction.\\
\textbf{Subcase 4.2} Let $L$-spectrum of $S$ be $\Big({\partial^L_{n-3} (G)-n},0\Big)$. For the $L$-spectrum of $\overline{G}$ and Fact \ref{fact}, we observe that the $L$-spectrum of $R$ and $S$ is same and is given by $\Big({(\partial^L_1 (G)-n)}^{(p-1)},0\Big)$. Clearly, $p\geq 3$, otherwise we have  $\partial^L_1 (G)= \partial^L_{n-3} (G)$, which is a contradiction. Hence, in this case $\overline{G}\cong K_P \cup K_{p} \cup K_2$ or  $G\cong K_{p,p,2}$, $p=\frac{n-2}{•2}\geq 3$.   \\
\textbf{Case 5.} Let $p=|S|\geq 3$. From the $L$-spectrum of $\overline{G}$ and $|Q|\geq |R|\geq |S|=p\geq 3$, we see that $\partial^L_{1} (G)-n $ is contained in the $L$-spectrum of all of $Q$, $R$ and $S$. Again, by Fact \ref{fact}, we easily get  $\overline{G}\cong K_p \cup K_{p} \cup K_{p-2}\vee \overline{K_2} $. Hence  $G\cong K_{p,p,p}+e$, $p=\frac{n}{•3}\geq 3$.\\
\indent For the converse, taking Lemmas \ref{L3}, \ref{L4}, \ref{L6} and \ref{L7} into consideration, it can be easily seen that the $D^L$-spectrum of $G\cong SK_{n,n-2}+e$, $G\cong K_{n,n-3}-e $, $G\cong K_{p,p,1}+e$, $G\cong K_{p,p,2}$ and $G\cong K_{p,p,p}+e$ are $\Big((2n-2)^{(n-4)},2n-4,n^{(2)},0\Big)$, $\Big((2n-3)^{(n-4)},n+2,n^{(2)},0\Big)$, $\Big((\frac{3n-1}{•2})^{(n-4)},\frac{3n-5}{•2},n^{(2)},0\Big)$, $\Big((\frac{3n-2}{•2})^{(n-4)},n+2,n^{(2)},0\Big)$ and $\Big((\frac{4n}{•3})^{(n-4)},\frac{4n-6}{•3},n^{(2)},0\Big)$, respectively.  \qed

\indent To completely characterize the connected graphs with $n\geq 5$ vertices, where distance Laplacian spectral radius has multiplicity $n-4$, we need to find the solution for the Cases (c) and (d), which are left as open problems.

\section{Multiplicity of any distance Laplacian eigenvalue}\label{3}

 To prove the next theorem, we need the following lemma.
 \begin{lemma}\label{LP2} {\em \cite{A4}}
 Let $G$ be a graph on $n\geq 3$ vertices whose distinct Laplacian eigenvalues are $0<\alpha <\beta$. The multiplicity of $\alpha$ is $n-2$ if and only if $G$ is one of the graphs $K_{\frac{n}{•2}, \frac{n}{•2}}$ or $S_n$.
\end{lemma}
\begin{theorem}\label{T2}
Let $G$ be a connected graph of order $n\geq 5$ having $\partial^L_1 (G)$ of multiplicity one. Then $\partial^L_{n-1} (G)=n$ with multiplicity 2 and $m(\partial^L_{n-3} (G))=n-4$ if and only if $G\cong S_3 \vee (K_2 \cup K_2)$ for $n=7$ or $G\cong (K_{n-1}\cup K_1 )+2e$ or $G\cong K_2 \vee (K_{\frac{n-2}{•2}} \cup K_{\frac{n-2}{•2}})$ or  $G\cong K_{\frac{n}{•4},\frac{n}{4•}}\vee (K_\frac{n}{•4} \cup K_\frac{n}{•4})$.
\end{theorem}
\noindent{\bf Proof.} Let $n$ be a distance Laplacian eigenvalue with multiplicity 2. Using Lemma \ref{L3}, we observe that $\overline{G}$ has 3 components and $diam(G)=2$. Let $\overline{G}\cong F\cup T\cup S $, where $|F|\geq |T|\geq |S|$. Using Lemmas \ref{L1} and \ref{L2}, we note that the $L$-spectrum of $\overline{G}$ is $\Big(\partial^L _1 (G)-n, (\partial^L _{n-3} (G)-n)^{(n-4)},0,0,0 \Big)$. We have the following cases to consider.\\
\textbf{ Case 1.} If $|S|=|T|=1$, then $L$-spectrum of $F$ is $\Big(\partial^L _1 (G)-n,(\partial^L _{n-3} (G)-n)^{(n-4)},0\Big)$. Applying Lemma \ref{LP2}, it is easy to see that either $F\cong S_{n-2}$ or $F\cong K_{\frac{n-2}{•2}, \frac{n-2}{•2}}$. Since $S$ and $T$ are both isomorphic to $K_1$, therefore $\overline{G} \cong S_{n-2}\cup K_1 \cup K_1$ or $\overline{G} \cong K_{\frac{n-2}{•2}, \frac{n-2}{•2}}\cup K_1 \cup K_1$ . This shows that $G $ is one of the graphs $ (K_{n-1}\cup K_1 )+2e$  and $K_2 \vee (K_{\frac{n-2}{•2}} \cup K_{\frac{n-2}{•2}})$.\\
\textbf{ Case 2.} If $|S|=1$ and $|T|=2$, then $L$-spectrum of $T$ is either $\Big(\partial^L _1 (G)-n,0 \Big)$ or $\Big(\partial^L _{n-3} (G)-n,0 \Big)$. The following two subcases arise.\\
\textbf{ Subcase 2.1.} Let the $L$-spectrum of $T$ be  $\Big(\partial^L _1 (G)-n,0 \Big)$. So  $T\cong K_2$ and the $L$-spectrum of  $F$ is $\Big((\partial^L _{n-3} (G)-n)^{(n-4)},0\Big)$. Using Fact \ref{fact}, from the above argument, we observe that $\partial^L _1 (G)-n=2$ and $\partial^L _{n-3} (G)-n=n-3\geq 2$.  Therefore, $ \partial^L _{n-3} (G)\geq \partial^L _1 (G)$, a contradiction.\\
\textbf{ Subcase 2.2.}  Let $L$-spectrum of $T$ be  $\Big(\partial^L _{n-3} (G)-n,0 \Big)$. Therefore, $L$-spectrum of $F$ is  $\Big(\partial^L _1 (G)-n,(\partial^L _{n-3} (G)-n)^{(n-5)},0\Big)$. If $n\geq 8$, then using Lemma \ref{LP2} and Fact \ref{fact}, we get $\partial^L _{n-3} (G)-n=2$ from $L$-spectrum of $T$ and $\partial^L _{n-3} (G)-n=1 $ or  $\partial^L _{n-3} (G)-n=\frac{n-3}{•2} $ from $L$-spectrum of $F$,  which is a contradiction. If $n=5$, we have $\partial^L _{n-3} (G)=\partial^L _1 (G)$, a contradiction. If $n=6$, using Lemma \ref{LP2},  $\partial^L _{n-3} (G)-n=1$ from $L$-spectrum of $F$ which is a again a contradiction. If $n=7$, using the same arguments as above, $F$ is one of the graphs $S_4$ or $K_{2,2}$. $F\cong S_4$ gives a contradiction while as $F\cong K_{2,2}$ shows that $G\cong S_3 \vee (K_2 \cup K_2)$.\\
\textbf{ Case 3.} If $|S|=1$ and $k=|T|\geq 3$, then from the $L$-spectrum of $\overline{G}$, we see that either $F$ contains three distinct Laplacian eigenvalues or $T$  contains three distinct Laplacian eigenvalues. It suffices to consider one of the two cases. Without loss of generality, assume that $T$ contains three distinct Laplacian eigenvalues. So the Laplacian spectrum of $F$ is $\Big((\partial^L _{n-3} (G)-n)^{(n-k-2)},0\Big)$ and  the Laplacian spectrum of $T$ is $\Big(\partial^L _1 (G)-n,(\partial^L _{n-3} (G)-n)^{(k-2)},0\Big)$. Applying Lemma \ref{LP2} and Fact \ref{fact}, we get  $\partial^L _{n-3} (G)-n=1$ or $ \frac{k}{•2}$ from $L$-spectrum of $T$ and $\partial^L _{n-3} (G)-n=n-k-1$ from $L$-spectrum of $F$, a contradiction. \\
\textbf{ Case 4.} If $|S|\geq 2$, then we observe from the $L$-spectrum of $\overline{G}$ that only one component among the $F,~T$ and $S$ contains $\partial^L _1 (G)-n$ as a Laplacian eigenvalue. The $L$-spectrum of the remaining two is the same. Note that for $n=6$, $\partial^L _{n-3} (G)=\partial^L _1 (G)$, which is a contradiction. Let $b$ be the cardinality of the component containing $\partial^L _1 (G)-n$ as Laplacian eigenvalue. For $n\geq 7$, if $b=2$, we observe that $\partial^L _1 (G)-n=2$ from the $L$-spectrum of the component containing $\partial^L _1 (G)-n$ as Laplacian eigenvalue and $\partial^L _{n-3} (G)-n\geq 2$ from the spectrum of  remaining two components, a contradiction. Let $b\geq 3$, using Lemma \ref{LP2}, $\partial^L _{n-3} (G)-n=1$ or $\partial^L _{n-3} (G)-n=\frac{b}{•2} $. $\partial^L _{n-3} (G)-n=1$ gives a contradiction and $\partial^L _{n-3} (G)-n=\frac{b}{•2} $ shows that $n=2b$ and $ \overline{G} \cong K_{\frac{b}{•2},\frac{b}{2•}}\cup K_\frac{b}{•2} \cup K_\frac{b}{•2}$ so that $G\cong K_{\frac{n}{•4},\frac{n}{4•}}\vee (K_\frac{n}{•4} \cup K_\frac{n}{•4})$.   \\
\indent Conversely, noting that all the graphs in the statement of the theorem are of diameter two and using Lemmas \ref{L1} and \ref{L2}, it is easy to see that the $D^L$-spectrum of $G\cong S_3 \vee (K_2 \cup K_2)$, $G\cong K_{n}-(n-3)e$, $G\cong K_2 \vee (K_{\frac{n-2}{•2}} \cup K_{\frac{n-2}{•2}})$ and $G\cong K_{\frac{n}{•4},\frac{n}{4•}}\vee (K_\frac{n}{•4} \cup K_\frac{n}{•4})$ are $\Big(11,9^{(3)},7^{(2)},0\Big)$, $\Big(2n-2,(n+1)^{(n-4)},n^{(2)},0\Big)$, $\Big(2n-2,(\frac{3n-2}{•2})^{(n-4)},n^{(2)},0\Big)$ and $\Big(\frac{3n}{•2},(\frac{5n}{•4})^{(n-4)},n^{(2)},0\Big)$, respectively. \qed

\indent Now, we will completely determine the graphs for which $n$ is a distance Laplacian eigenvalue of multiplicity $n-4$.

\begin{theorem}\label{T3}
Let $G$ be a connected graph with order $n\geq 5 $. Then\\
$ (a)$  $ m(\partial^L_1 (G))=3$ and  $ \partial^L_{n-1} (G)=n$ with multiplicity $n-4$ if and only if $G\cong K_{4,1,\dots,1}$ or $G\cong K_{2,2,2,1,\dots,1}$.\\
$ (b)$  $ m(\partial^L_1 (G))=2$ and  $ \partial^L_{n-1} (G)=n$ with multiplicity $n-4$ if and only if $G\cong SK_{n,4}+e$ or $G\cong K_{3,2,1,\dots,1}$.\\
$ (c)$  $ m(\partial^L_1 (G))=1$ and  $ \partial^L_{n-1} (G)=n$ with multiplicity $n-4$ if and only if $G$ is isomorphic to any one of the following graphs, $ \overline{S_4 \cup (n-4)K_1}$ or    $\overline{C_4 \cup (n-4)K_1}$ or  $ \overline{P_4 \cup (n-4)K_1}$ or $ \overline{Ki_{4,3}  \cup (n-4)K_1}$ or  $\overline{ S_3 \cup K_2 \cup (n-5)K_1}$.
\end{theorem}
\noindent{\bf Proof.} {\bf (a).} Let $n$ be a distance Laplacian eigenvalue of $G$ with multiplicity $n-4$ and  $ m(\partial^L_1 (G))=3$. Using Lemma \ref{L3}, $\overline{G}$ is disconnected with $n-3$ components and $diam(G)=2$. Applying Lemmas \ref{L1} and  \ref{L2}, the $L$-spectrum  of $\overline{G}$ is $\Big((\partial^L _1 (G)-n)^3 ,0,0,\dots,0\Big)$. From the $L$-spectrum of $\overline{G}$, we observe that $\overline{G}$ has exactly one non-zero Laplacian eigenvalue. So all the components of $\overline{G}$ are either isolated vertices or complete graphs of same order. Therefore, $\overline{G} \cong K_4 \cup (n-4)K_1$ or  $\overline{G} \cong K_2 \cup K_2 \cup K_2 \cup (n-6)K_1$. This further implies that $G\cong K_{4,1,\dots,1}$ or $G\cong K_{2,2,2,1,\dots,1}$.\\
\indent Conversely, by using Lemma \ref{L4}, we see that the $D^L$-spectrum of  $ SK_{n,4}$ and $ K_{2,2,2,1,\dots,1}$ are respectively, $\Big((n+4)^3 ,n^{(n-4)},0\Big)$ and $\Big((n+2)^3 ,n^{(n-4)},0\Big)$. \\
\noindent {\bf (b).} Now, let $n$ be a distance Laplacian eigenvalue of $G$ with multiplicity $n-4$ and  $ m(\partial^L_1 (G))=2$. So by using the same argument as in (a), we observe that $\overline{G}$ is disconnected with $n-3$ components, $diam(G)=2$  and the $L$-spectrum  of $\overline{G}$ is $\Big((\partial^L _1 (G)-n)^{2} ,\partial^L _3 (G)-n,0,\dots,0\Big)$. Let $\overline{G}\cong G_1 \cup G_2 \cup \dots \cup G_{n-3}$, where $G_i$, $i=1,2,\dots,n-3$, are the components of $\overline{G}$. Clearly, either one or two or at most three components of $\overline{G}$ can contain all the non-zero Laplacian eigenvalues. So we have the following cases to consider.\\
\textbf{Case 1. } Only one component contains all the non-zero Laplacian eigenvalues of $\overline{G}$. Without loss of generality, let $G_1$ contain  all the non-zero Laplacian eigenvalues of $\overline{G}$. So the $L$-spectrum of $G_1$ is  $\Big((\partial^L _1 (G)-n)^{2} ,\partial^L _3 (G)-n,0\Big)$. Note that there are only six connected graphs of order 4 as shown in Figure 1. By Fact \ref{fact}, only one graph $(K_4 -e)$ has Laplacian spectral radius of multiplicity 2. Hence in this case $ \overline{G} \cong( K_4 -e)\cup (n-4)K_1$, so that $G\cong SK_{n,4}+e$.\\
\textbf{Case 2. } Now, let two components contain all the non-zero Laplacian eigenvalues of $\overline{G}$. Without loss of generality, let $G_1$ and $G_2$ contain  all the non-zero Laplacian eigenvalues of $\overline{G}$. Also, assume that  $G_1$ contains two out of three non-zero Laplacian eigenvalues and the remaining one is contained in  $G_2$. Now, consider the following subcases.\\
\textbf{Subase 2.1. } First, let the $L$-spectrum of $G_1$ be  $\Big((\partial^L _1 (G)-n)^{2} ,0\Big)$ and $L$-spectrum of $G_2$ be  $\Big(\partial^L _3 (G)-n ,0\Big)$. Therefore, $G_1 \cong K_3$ and $G_2 \cong K_2$ which further implies that
$ \overline{G} \cong K_3\cup K_2\cup (n-5)K_1$. Hence, $G\cong K_{3,2,1,\dots,1}$ in this case.\\
 \textbf{Subase 2.2. } Let the $L$-spectrum of $G_1$ be  $\Big(\partial^L _1 (G)-n ,\partial^L _3 (G)-n,0\Big)$ and $L$-spectrum of $G_2$ be  $\Big(\partial^L _1 (G)-n ,0\Big)$. Using Lemma \ref{LP2} and Fact \ref{fact},  from $L$-spectrum of $G_1$, we get  $\partial^L _1 (G)-n=3$, and from $L$-spectrum of $G_2$ we get $\partial^L _1 (G)-n=2$. Clearly, this is a contradiction.\\
 \textbf{Case 3. } Let three components contain all the non-zero Laplacian eigenvalues of $\overline{G}$. Suppose that  $G_1$, $G_2$  and $G_3$ be those components. Without loss of generality, let their spectrum be  $\Big(\partial^L _1 (G)-n ,0\Big)$,  $\Big(\partial^L _1 (G)-n ,0\Big)$ and  $\Big(\partial^L _3 (G)-n ,0\Big)$. Using Fact \ref{fact}, we get  $\partial^L _1 (G)=\partial^L _3 (G)$, which is a contradiction.\\
 \indent Conversely, using Lemmas \ref{L6} and \ref{L7} and the fact that the  trace of a matrix equals to the sum of all its eigenvalues, we see that that  $D^L$-spectrum  of $ SK_{n,4}+e$ is $\Big((n+4)^2 ,n+2,n^{(n-4)},0\Big)$. From Lemma \ref{L4}, $D^L$-spectrum  of   $K_{3,2,1,\dots ,1}$ is $\Big((n+3)^2 ,n+2,n^{(n-4)},0\Big)$.\\
 \noindent {\bf (c).} Using the same arguments as in part (a) and (b), we observe that the $L$-spectrum  of $\overline{G}$ is  $\Big(\partial^L _1 (G)-n  ,\partial^L _2 (G)-n,\partial^L _3 (G)-n,0,\dots,0\Big)$ with $n-3$ components and $diam(G)=2$. We have the following two cases to consider.\\
 \textbf{Case 4. } Let $\partial^L _2 (G)=\partial^L _3 (G)$. Now, if only one component, say $G_1$, of  $\overline{G}$ contains all its non-zero Laplacian eigenvalues, then clearly $L$-spectrum  of $G_1$ is given by  $ \Big(\partial^L _1 (G)-n  ,(\partial^L _2 (G)-n)^(2) ,0\Big)$. Among all the six connected graphs on four vertices as shown in Figure 1, only the star and the cycle have second largest Laplacian eigenvalue of multiplicity two. Thus, either $G_1 \cong S_4$ or  $G_1 \cong C_4$, so that either $\overline{G}\cong S_4 \cup (n-4)K_1$ or $\overline{G}\cong C_4 \cup (n-4)K_1$. Therefore, $G\cong \overline{S_4 \cup (n-4)K_1}$ or    $G\cong \overline{C_4 \cup (n-4)K_1}$. Finally, if two or three components of  $\overline{G}$ contain all its non-zero Laplacian eigenvalues, then proceeding similarly as in (b), we arrive at a contradiction in both the cases.\\
 \textbf{Case 5. } Let $\partial^L _2 (G) \neq  \partial^L _3 (G)$. If only one component, say $G_1$ of  $\overline{G}$ contains all its non-zero Laplacian eigenvalues, then clearly $L$-spectrum  of $G_1$ is given by $ \Big(\partial^L _1 (G)-n  ,\partial^L _2 (G)-n ,\partial^L _3 (G)-n ,0\Big)$. Among all the six connected graphs on four vertices as shown in Figure 1, only the path $P_4$ and the kite $Ki_{4,3}$ have all the three non-zero Laplacian eigenvalues different. Thus, $G_1 \cong P_4$ or $G_1 \cong Ki_{4,3} $, so that $\overline{G}\cong P_4 \cup (n-4)K_1$ or $\overline{G}\cong Ki_{4,3} \cup (n-4)K_1$. Therefore, $G\cong \overline{P_4 \cup (n-4)K_1}$ or $G\cong \overline{Ki_{4,3}  \cup (n-4)K_1}$. Using the arguments as in part (b), the only case that remains to be discussed is when two components, say $G_1$ and $G_2$,  of  $\overline{G}$ contain all its non-zero Laplacian eigenvalues and $L$-spectrum of one, say  $G_1$, is $ \Big(\partial^L _1 (G)-n  ,\partial^L _3 (G)-n ,0\Big)$ and of $G_2$ is $ \Big(\partial^L _2 (G)-n  ,0\Big)$. Clearly, $G_1 \cong S_3$ and $G_2 \cong K_2$, so that  $\overline{G}\cong S_3 \cup K_2 \cup (n-5)K_1$. Therefore, $G\cong \overline{ S_3 \cup K_2 \cup (n-5)K_1}$.\\\\
\begin{figure}
\centering
\includegraphics[scale=.2]{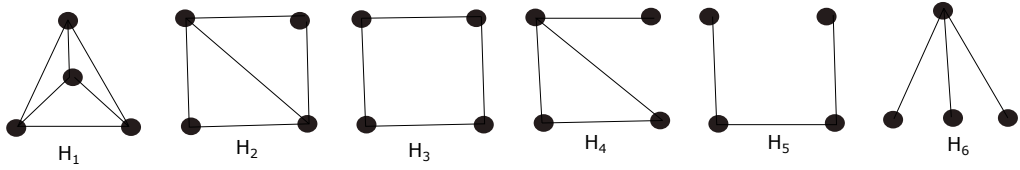}
\caption{All connected graphs on four vertices}
\end{figure}
\indent Conversely, note that the Laplacian spectrum of graphs $ S_4 \cup (n-4)K_1$,    $ C_4 \cup (n-4)K_1$,  $ P_4 \cup (n-4)K_1$, $ Ki_{4,3}  \cup (n-4)K_1$ and  $ S_3 \cup K_2 \cup (n-5)K_1$ are $\Big(4 ,1^{(2)},0^{(n-3)}\Big)$,  $\Big(4sin^2 (\frac{\pi}{4}), 4sin^2 (\frac{2\pi}{4}),4sin^2 (\frac{3\pi}{4}) ,0^{(n-3)}\Big)$, $\Big(4sin^2 (\frac{\pi}{8}) ,4sin^2 (\frac{2\pi}{8}),4sin^2 (\frac{3\pi}{8}),0^{(n-3)}\Big)$, $\Big(4 ,3,1,0^{(n-3)}\Big)$ and $\Big(3 ,2,1,0^{(n-3)}\Big)$, respectively. Also the complement of all these graphs are of diameter two. Using Lemmas \ref{L1} and \ref{L2}, we see that the $D^L$-spectrum of
 $ \overline{S_4 \cup (n-4)K_1}$,    $ \overline{C_4 \cup (n-4)K_1}$,  $ \overline{P_4 \cup (n-4)K_1}$, $ \overline{Ki_{4,3}  \cup (n-4)K_1}$ and  $\overline{ S_3 \cup K_2 \cup (n-5)K_1}$ are $\Big(n+4, (n+1)^{2},n^{(n-4)},0\Big)$, $\Big(n+4, (n+2)^2,n^{(n-4)},0\Big)$, $\Big(n+4sin^2 (\frac{3\pi}{8}) ,n+4sin^2 (\frac{2\pi}{8}),n+4sin^2 (\frac{\pi}{8}),n^{(n-4)},0\Big)$,  $\Big(n+4 ,n+3,n+1,n^{(n-4)},0\Big)$  and  $\Big(n+3 ,n+2,n+1,n^{(n-4)},0\Big)$, respectively. \qed

\noindent{\bf Acknowledgements.}  The research of S. Pirzada is supported by SERB-DST, New Delhi under the research project number CRG/2020/000109. The research of Saleem Khan is supported by MANUU.

\noindent{\bf Data availability} Data sharing is not applicable to this article as no data sets were generated or analyzed
during the current study.

\end{document}